\documentclass[11pt,reqno]{amsart}
\usepackage{color}

\usepackage{amsmath}
\usepackage{amsfonts,amscd}
\usepackage{amssymb}
\usepackage{url}
\usepackage{hyperref}

\usepackage[english]{babel}
\usepackage{booktabs}
\usepackage{tikz}

\theoremstyle{plain}
\newtheorem{theorem}                 {Theorem}      [section]

\newtheorem{proposition}  [theorem]  {Proposition}

\theoremstyle{definition}

\newtheorem{definition}   [theorem]  {Definition}

\newtheorem{example}      [theorem]  {Example}

\numberwithin{equation}{section}

\def \cn{{\mathbb C}}
\def \C{{\mathbb C}}

\def \H{{\mathbb H}}

\def \rn{{\mathbb R}}

\def \zn{{\mathbb Z}}

\def \C{\mathcal C}
\def \E{\mathcal E}

\def \H{\mathcal H}

\def \Z{\mathcal Z}

\def\nab#1#2{\hbox{$\nabla$\kern -.3em\lower 1.0 ex
		\hbox{$#1$}\kern -.1 em {$#2$}}}
\def\hatnab#1#2{\hbox{$\nabla$\kern -.3em\lower 1.0 ex
		\hbox{$#1$}\kern -.1 em {$#2$}}}

\def \ip #1#2{\langle #1,#2 \rangle}

\def \SL2{\widetilde{\text{\bf SL}}_{2}(\rn)}

\DeclareMathOperator{\grad}{grad}

\numberwithin{equation}{section}
\allowdisplaybreaks

\begin{document}

	

\author{Sigmundur Gudmundsson}
\address{Mathematics, Faculty of Science\\
Lund University\\
Box 118, Lund 221 00\\
Sweden}
\email{Sigmundur.Gudmundsson@math.lu.se}

\author{Leonard Nygren Löhndorf}
\address{Mathematics, Faculty of Science\\
Lund University\\
Box 118, Lund 221 00\\
Sweden}
\email{leonard.nygren.lohndorf@gmail.com}

\title
[New Minimal Submanifolds]
{New Minimal Surfaces of the Sphere $S^4$\\ and the Hyperbolic Space $H^4$\\  via Harmonic Morphisms}

\begin{abstract}
In this work we introduce a new method for the construction of minimal submanifolds of codimension two in even dimensional spheres and hyperbolic spaces.  This is based on the theory of complex-valued harmonic morphisms.  This gives the first explicit examples of such maps defined on the sphere $S^4$ and the hyperbolic space $H^4$.
\end{abstract}
	

\maketitle

\section{Introduction}
\label{section-introduction}

The classical Hopf map $\phi:S^3\to S^2$ is named after the distinguished mathematician Heinz Hopf (1894-1971).  This is a smooth map between the compact three dimensional unit sphere $S^3$ in $\cn^2\cong\rn^4$ and the compact two dimensional $S^2$ in $\rn^3\cong\cn\times\rn$.  In the complex coordinates $(z,w)\in\cn^2$ it is given explicitly by 
$$
\phi:(z,w)\mapsto (2z\bar w,|z|^2-|w|^2).
$$
This has later been generalised in several different directions and been important  both in  topology and differential geometry.
\smallskip

For the theory of harmonic morphisms this interesting map has played an important role. Early on, Baird and Eells prove in \cite{Bai-Eel} that the Hopf map is a harmonic morphism.  These are known to be rather rare between compact Riemannian manifolds.  Identifying $S^2$ with the Riemann sphere $\hat\cn=\cn\cup\{\infty\}$ the above presentation of the Hopf map $\phi:S^3\to\hat\cn$ takes the form
$$\phi:(z,w)\mapsto \frac zw=\frac{x_1+ ix_2}{x_3+i x_4},$$
where $z=x_1+ i\, x_2$, $w=x_3+ i\, x_4$ and $(x_1,x_2,x_3,x_4)\in\rn^4\cong\cn^2$.
\smallskip 

Later the author of \cite{Gud-7} employed the well-known duality between the three dimensional sphere $S^3$ and the three dimensional hyperbolic space $H^3$, as Riemannian symmetric spaces, to construct the globally defined harmonic morphism $\phi^*:H^3\to\cn$ satisfying 
$$\phi^*:(x_1,x_2,x_3,x_4)\mapsto \frac{x_1+ ix_2}{x_3-x_4}.$$
Here the hyperbolic space $H^3$ is the standard model embedded into the classical Minkowski space $\rn^4_1$.  This was then easily extended to the odd dimensional cases of $S^{2n+1}$ and $H^{2n+1}$ by noting that $S^{2n+1}$ is the unit sphere on $\cn^{n+1}$,  allowing for basic applications of complex analysis.
\smallskip   

In his work \cite{Sve-1} Svensson gives a method to extend to the even dimensional cases of $H^{2n}$ and hence $S^{2n}$, for $n>2$. To our knowledge no explicit examples have so far been constructed in the cases of $S^4$ and $H^4$.  The main aim of this work is to close this gap in the current theory.  For the construction of our new examples see Sections \ref{section-minimal-S4} and \ref{section-minimal-H4}.

\section{Preliminaries}

Let $M$ and $N$ be two manifolds of dimensions $m$ and $n$,
respectively. Then a semi-Riemannian metric $g$ on $M$ gives rise
to the notion of a Laplace-Beltrami operator (alt. tension field) on $(M,g)$ and real-valued harmonic
functions $f:(M,g)\to\rn$. This can be generalised to the concept
of a harmonic map $\phi:(M,g)\to (N,h)$ between semi-Riemannian
manifolds being a solution to a non-linear system of partial
differential equations. 

\begin{definition}
A map $\phi:(M,g)\to (N,h)$ between semi-Riemannian manifolds is called a {\it harmonic morphism} if, for any harmonic function $f:U\to\rn$ defined on an open subset $U$ of $N$ with $\phi^{-1}(U)$ non-empty, the composition $f\circ\phi:\phi^{-1}(U)\to\rn$ is a harmonic function.
\end{definition}

For the theory of harmonic morphisms we refer to the excellent text \cite{Bai-Woo-book} and the regularly updated online bibliography \cite{Gud-bib}. The following characterization of harmonic morphisms between semi-Riemannian manifolds is due to Fuglede, see \cite{Fug-2}. For the definition of horizontal conformality we
refer to \cite{Bai-Woo-book}.

\begin{theorem}\cite{Fug-2}
A map $\phi:(M,g)\to (N,h)$ between semi-Riemannian manifolds is a harmonic morphism if and only if it is a horizontally (weakly) conformal harmonic map.
\end{theorem}

The next result generalises the corresponding well-known theorem
of Baird and Eells in the Riemannian case, see \cite{Bai-Eel}. It
gives the theory of harmonic morphisms a strong geometric flavour
and shows that the case when $n=2$ is particularly interesting. In
that case the conditions characterising harmonic morphisms are
independent of conformal changes of the metric on the surface
$N^2$.  For the definition of horizontal homothety we refer to
\cite{Bai-Woo-book}.

\begin{theorem}\cite{Gud-8}\label{theorem-minimal}
Let $(M^m,g)$ be a semi-Riemannian manifold, $(N^n,h)$ be Riemannian and $\pi:M\to N$ be a horizontally conformal submersion.  If
\begin{enumerate}
\item[i.] $n=2$, then $\pi$ is harmonic if and only if $\pi$ has minimal fibres,
\item[ii.] $n\ge 3$, then two of the following conditions imply the other,
\begin{enumerate}
\item $\pi$ is a harmonic map,
\item $\pi$ has minimal fibres, 
\item $\pi$ is horizontally homothetic.	
\end{enumerate}
\end{enumerate}
\end{theorem}

In what follows we are mainly interested in complex-valued functions $\phi,\psi:(M,g)\to\cn$ from semi-Riemannian manifolds. In this situation the metric $g$ induces the complex-valued {\it tension field} $\tau(\phi)$ and the gradient $\grad(\phi)$ with
values in the complexified tangent bundle $T^{\cn}M$ of $M$.  We extend the metric $g$ to be complex bilinear on $T^{\cn} M$ and define the symmetric bilinear {\it conformality operator} $\kappa$ by
$$\kappa(\phi,\psi)= g(\grad(\phi),\grad(\psi)).$$   
The harmonicity and horizontal conformality of $\phi:(M,g)\to\cn$ are given by the following relations
$$\tau(\phi)=0\ \ \text{and}\ \ \kappa(\phi,\phi)=0.$$

\begin{definition}\cite{Gud-Sak-1}\label{definition-eigen}
Let $(M,g)$ be a semi-Riemannian manifold.  Then a set
$$\E =\{\phi_i:M\to\cn\ |\ i\in I\}$$ 
of complex-valued functions is said to be an {\it eigenfamily} on $M$ if there exist complex numbers $\lambda,\mu\in\cn$ such that
$$\tau(\phi)=\lambda\cdot\phi\ \ \text{and}\ \ \kappa (\phi,\psi)=\mu\cdot\phi\,\psi,$$
for all $\phi,\psi\in\E$. 
\end{definition}

The next result shows that an eigenfamily on a semi-Riemannian manifold can be used to produce a variety of local harmonic morphisms.

\begin{theorem}\cite{Gud-Sak-1}\label{theorem-rational}
Let $(M,g)$ be a semi-Riemannian manifold and 
$$\E =\{\phi_1,\dots,\phi_n\}$$ 
be a finite eigenfamily of complex-valued functions on $M$. If $P,Q:\cn^n\to\cn$ are linearly independent homogeneous polynomials of the same positive degree then the quotient
$$\frac{P(\phi_1,\dots ,\phi_n)}{Q(\phi_1,\dots ,\phi_n)}$$ 
is a non-constant harmonic morphism on the open and dense subset
$$\Omega (Q)=\{p\in M\,|\, Q(\phi_1(p),\dots ,\phi_n(p))\neq 0\}.$$
\end{theorem}

\section{Harmonic Morphisms on the Euclidean $\rn^{n}$}

The study of the spherical harmonics is a beautiful part of our classical education in mathematics. There we study the complex-valued polynomials $P:\rn^{n}\to\cn$ which are harmonic i.e. those that satisfy the classical Laplace equation
$$\tau(P)=\sum_{k=1}^{n}\frac{\partial^2 P}{\partial x_k^2}=0.$$ 
These polynomials form a linear space $\H_{n}$ which has a natural direct linear decomposition 
$$\H_{n}=\bigoplus_{d=1}^\infty  \H^d_{n},$$
where $\H^d_{n}$ is the finite dimensional subspace of $\H_{n}$ consisting of the harmonic polynomials of degree $d\in\zn^+$.  

If $\phi,\psi:\rn^{n}\to\cn$ are complex-valued $C^1$-functions then the conformality operator $\kappa$ satisfies
$$\kappa(\phi,\psi)=\sum_{k=1}^{n}\frac{\partial\phi}{\partial x_k}\frac{\partial\psi}{\partial x_k}.$$ 
For this we now define the subset $\C\H_{n}^d$ of $\H_{n}^d$ by $$\C\H_{n}^d=\{P\in\H_{n}^d\,|\,\kappa(P,P)=0\}$$
i.e. the space of harmonic polynomials $P\in\H_{n}^d$ which are horizontally conformal.  Then $\C\H_{n}^d$ is a complex cone in $\H_{n}^d$ i.e. if $P\in\C\H_{n}^d$ then so is $\alpha\cdot P$, for all non-zero complex numbers $\alpha\in\cn^*$.

\begin{example}
Let us now take a special look at the linear space $\H_{n}^1$ with the canonical basis 
$$\{x_1,x_2,\dots ,x_{n}\}.$$
Then every complex-valued element $P\in\H_{n}^1$ can be expressed as
$$P(x)=\sum_{k=1}^{n}c_k\, x_k,$$
with complex coefficients $c_1,c_2,\dots ,c_{n}\in\cn$. A simple calculation shows that for the conformality operator $\kappa$ we have $$\kappa(P,P)=c_1^2+c_2^2+\cdots +c_{n}^2.$$
Then it is clear that the above harmonic element $P$ belongs to $\C\H_n^1$ if and only if 
$c_1^2+c_2^2+\cdots +c_{n}^2=0$ i.e. the complex vector $(c_1,c_2,\dots,c_{n})\in\cn$ is isotropic.
\end{example}
\medskip

\begin{example}
The linear space $\H^2_2$ of second order harmonic polynomials on the Euclidean space $\rn^2$ is $2$-dimensional.  For this we have the following natural  basis
$$\{x_1^2-x_2^2,x_1x_2\}.$$
This means that every complex-valued harmonic element  $P\in\H^2_2$ can be expressed by
$$P(x)=a\,(x_1^2-x_2^2)+b\,(x_1x_2),$$
where the coefficients $a,b\in\cn$ are complex numbers.
An elementary calculation shows that every such element $P:\rn^2\to\cn$ is horizontally conformal i.e. $\kappa (P,P)=0$ if and only if
$$a^2+4\cdot b^2=0.$$
This implies that such a function $P$ is of the form $$P(x_1,x_2)=a\cdot(x_1\pm i\cdot x_2)^2$$ i.e. it is either  holomorphic or anti-holomorphic on $\cn\cong\rn^2$.
\end{example}
\medskip

\begin{example}\label{example-p}
The linear space $\H^2_3$ of second order harmonic polynomials on the Euclidean space $\rn^3$ is $5$-dimensional.  For this we have the following natural  basis
$$\{x_1^2-x_2^2,x_2^2-x_3^2,x_1x_2,x_1x_3,x_2x_3\}.$$
This means that every complex-valued harmonic element  $p\in\H^2_3$ can be expressed as
$$p(x)=a_{1}(x_1^2-x_2^2)+a_{2}(x_2^2-x_3^2)+b_{1}(x_1x_2)+b_{2}(x_1x_3)+b_{3}(x_2x_3),$$
where the coefficients $a_{1},a_{2},b_{1},b_{2},b_{3}\in\cn$ are complex numbers.
An elementary calculation shows that $p:\rn^3\to\cn$ is horizontally conformal i.e. $\kappa (p,p)=0$ if and only if
$$a_{1}^2 + b_{1}^2 + b_{2}^2=0,\ a_{1}a_{2} = -b_{2}^2\ \ \text{and}\ \ 
a_{1}b_{3} = b_{1} b_{2}.$$

Note that these polynomials form a set parametrised by an algebraic variety of complex dimension two.
\end{example}

\begin{example}\label{example-PQ}
For a positive integer $n\in\zn^+$ let $\cn^n$ denote the standard Euclidean $n$-dimensional complex vector space of real dimension $2n$.  Let $P_d,Q_d:\cn^n\to\cn$ be two linearly independent homogeneous polynomials in the complex variables $(z_1,z_2,\dots,z_n)\in\cn^n$ of the same positive degree $d\in\zn^+$.  Then it is well-known that they are harmonic i.e. $\tau(P_d)=0$ and $\tau(Q_d)=0$ and furthermore we have 
$$\kappa(P_d,P_d)=0,\ \kappa(P_d,Q_d)=0,\ \kappa(Q_d,Q_d)=0.$$
Note that in this case the tension field $\tau$ and the conformality operator $\kappa$ on $\cn^n$ are given by
$$\tau(\phi)=4\cdot\sum_{k=1}^n\frac{\partial^2 \phi}{\partial z_k\partial\bar z_k}
\ \ \text{and}\ \ \kappa(\phi,\psi)=2\cdot\sum_{k=1}^n\big(\frac{\partial\phi}{\partial z_k}\frac{\partial\psi}{\partial\bar z_k}+\frac{\partial\phi}{\partial \bar z_k}\frac{\partial\psi}{\partial z_k}\big).$$
\end{example}

\begin{example}
Let $P_d,Q_d:\cn^n\to\cn$ be the polynomials defined in Example \ref{example-PQ} and $p:\rn^3\to\cn$ as in Example \ref{example-p} with 
$$a_{1}^2 + b_{1}^2 + b_{2}^2=0,\ a_{1}a_{2} = -b_{2}^2\ \ \text{and}\ \ 
a_{1}b_{3} = b_{1} b_{2}.$$
For an {\it even} positive integer $d$ we define the complex-valued functions $\hat\Phi_d:(\rn^3\times\cn^n)\setminus\Z(Q_d)\to\cn$ with
$$\hat\Phi_d:(x_1,x_2,x_3,z_1,z_2,\dots ,z_n)\mapsto \frac {p(x_1,x_2,x_3)^{d/2}+P_d(z_1,z_2,\dots,z_n)}{Q_d(z_1,z_2,\dots,z_n)}.$$
Here $\Z(Q_d)$ is the zero set of the polynomial $Q_d$ i.e. 
$$\Z(Q_d)=\{(x,z)\in\rn^3\times\cn^n\,|\,Q_d(z_1,z_2,\dots,z_n)=0\}.$$
For an {\it odd} positive integer $d$ we define the complex-valued function by $\hat\Phi_d:(\rn^3\times\cn^n)\setminus(\Z(Q_d)\cup\rn^-_0)\to\cn$ with the same formula as in the even case.  By $p^{1/2}=\sqrt p$ we mean the standard square root well defined on the open subset $\cn\setminus\{z\in\cn\cap\rn\,|\, z\le 0\}$ of $\cn$. 

For the product space $\rn^3\times\cn$ the tension field $\tau$ and the conformality operator $\kappa$ are given by
$$\tau(\phi)=\sum_{k=1}^3\frac{\partial^2 \phi}{\partial x_k^2}+4\cdot\sum_{k=1}^n\frac{\partial^2 \phi}{\partial z_k\partial\bar z_k}$$
and 
$$\kappa(\phi,\psi)=\sum_{k=1}^3\frac{\partial\phi}{\partial x_k}\frac{\partial\psi}{\partial x_k}
+2\cdot\sum_{k=1}^n\big(\frac{\partial\phi}{\partial z_k}\frac{\partial\psi}{\partial\bar z_k}+\frac{\partial\phi}{\partial \bar z_k}\frac{\partial\psi}{\partial z_k}\big).$$

With this at hand it is easily checked that if $P=p^{d/2}+P_d$ and $Q=Q_d$ then 
$$\tau(P)=0,\ \tau(Q)=0,\ \kappa(P,P)=0,\ \kappa(P,Q),\ \kappa(Q,Q)=0.$$ 
This means that $P,Q$ belong to the same eigenfamily on their open domains in $\rn^3\times\cn^n\cong\rn^{2n+3}$.  According to Theorem \ref{theorem-rational} the maps $\hat\Phi_d$ are all harmonic morphisms.
\end{example}

\section{Harmonic Morphisms on the Spheres $S^{2n+2}$}\label{section-hm-S2n+2}

Let $S^{2n+2}$ be the standard round sphere in $\rn^{2n+3}$ and $\pi:\rn^{2n+3}\setminus\{0\}\to S^{2n+2}$ be the radial projection given by $\pi:x\mapsto x/|x|$.  Then $\pi$ is a well-known submersive  harmonic morphism.

The next useful result is a special case of Proposition 1 in \cite{Gud-8}.

\begin{proposition}\label{proposition-invariance}
Let $\pi:(M,g)\to(N,h)$ be a non-constant harmonic morphism between two Riemannian manifolds. Further let $\phi:(N,h)\to\cn$ be a complex-valued function and $\hat \phi:(M,g)\to\cn$ be the composition $\hat\phi=\phi\circ\pi$.  Then $\phi$ is a harmonic morphism if and only if $\hat\phi$ is a harmonic morphism.
\end{proposition}

The complex-valued functions $\hat\Phi_d:\hat U_d\subset M\to\cn$ are invariant under the action of $\rn^+$ on $\rn^{2n+3}\setminus\{0\}$ i.e. $\hat\Phi_d(\lambda\cdot(x,z))=\hat\Phi_d((x,z))$ for all positive $\lambda\in\rn^+$.  Here $\hat U_d$ is the domain of $\hat\Phi_d$ in each case.
Hence they induce local complex-valued functions $\Phi_d:U_d=\pi(\hat U_d)\subset S^{2n+2}\to\cn$ on the corresponding round sphere $S^{2n+2}$.  According to Proposition \ref{proposition-invariance} the functions $\Phi_d$ are all local harmonic morphisms on the corresponding even-dimensional round sphere $S^{2n+2}$.

\section{Compact Minimal Surfaces in the Sphere $S^{4}$}\label{section-minimal-S4}
Let $a_1,a_2,b_1,b_2,b_3\in\cn$ be complex numbers such that 
$$a_{1}^2 + b_{1}^2 + b_{2}^2=0,\ a_{1}a_{2} = -b_{2}^2,\ a_{1}b_{3} = b_{1} b_{2}\ \ \text{and}\ \ a_1b_1b_2\neq 0.$$
Then define the polynomial functions $P,Q:(\rn^3\times\cn)\setminus\{0\}$ by 
$$P(x)=a_{1}(x_1^2-x_2^2)+a_{2}(x_2^2-x_3^2)+b_{1}(x_1x_2)+b_{2}(x_1x_3)+b_{3}(x_2x_3)$$
and 
$$Q(z)=z^2=(x_4+i\cdot x_5)^2.$$
According to our investigations above we know that the quotient $\hat\Phi=P/Q$ is a harmonic morphism on the open and dense subset $\rn^5\setminus(\Z(Q)\cup \{0\})$ of $\rn^5$.  The complex gradient $\grad(\hat\Phi)$
is given by 
$$\grad(\hat\Phi)=(
\frac{\partial\hat\Phi}{\partial x_1},
\frac{\partial\hat\Phi}{\partial x_2},
\frac{\partial\hat\Phi}{\partial x_3},
\frac{\partial\hat\Phi}{\partial x_4},
\frac{\partial\hat\Phi}{\partial x_5}),$$
where 
{\small $$\frac{\partial\hat\Phi}{\partial x_1}=
\frac{2a_1x_1 + b_1x_2 + b_2x_3}{(x_4 + i x_5)^2},$$
$$\frac{\partial\hat\Phi}{\partial x_2}=
\frac{b_1(a_1x_1+2b_1x_2 + b_2x_3)}{a_1(x_4 + i x_5)^2},$$
$$\frac{\partial\hat\Phi}{\partial x_3}=
\frac{b_2(a_1x_1 + b_1x_2 + 2b_2x_3)}{a_1(x_4 +i x_5)^2},$$
$$\frac{\partial\hat\Phi}{\partial x_4}=
\frac{-2a_1^2(x_1^2 -x_2^2) + 2b_2^2(x_2^2 - x_3^2) - 2a_1b_1x_1x_2 - 2a_1b_2x_1x_3 - 2b_1b_2x_2x_3}{a_1(x_4 +i x_5)^3},$$
$$\frac{\partial\hat\Phi}{\partial x_5}=-i\cdot\frac{\partial\hat\Phi}{\partial x_4}.$$}

Note that for the following determinant we have
$$
\det 
\begin{bmatrix}
2a_1 &  b_1 &  b_2 \\
 a_1 & 2b_1 &  b_2 \\
 a_1 &  b_1 & 2b_2
\end{bmatrix}
=4a_1b_1b_2\neq 0.
$$
This tells us that the only solutions to the equation $\grad(\hat\Phi)=0$ are the elements of the punctured $(x_4,x_5)$-plane i.e. $$\Pi=\{(0,0,0,x_4,x_5)\,|\,x_4^2+x_5^2\neq 0\}.$$
For a point $x\in\Pi$ we have $\hat\Phi(x)=0$.  This shows that the only critical value for $\hat\Phi$ is $0\in\cn$.  This means that for any non-zero complex number $\alpha\in\cn^*$ in the image of $\hat\Phi$ the fibre $\hat\Phi^{-1}(\{\alpha\})$ is a regular minimal submanifold of $\rn^5$ of codimension two.  Define the complex-valued function $F:\rn^5\to\cn$ by 
$$F(x_1,x_2,x_3,x_4)=P(x_1,x_2,x_3)-\alpha\cdot Q(x_4,x_5),$$
then 
$$\frac{\partial F}{\partial x_1}=a_1(2a_1x_1 +  b_1x_2 +  b_2x_3),$$
$$\frac{\partial F}{\partial x_2}=b_1( a_1x_1 + 2b_1x_2 +  b_2x_3),$$
$$\frac{\partial F}{\partial x_3}=b_2( a_1x_1 +  b_1x_2 + 2b_2x_3),$$
$$\frac{\partial F}{\partial x_4}=2\alpha (x_4 + i\,x_5)\ \ \text{and}\ \ \frac{\partial F}{\partial x_5}=-2i\alpha (x_4 + i\,x_5).$$
Then the complex gradient $\grad F(x)=0$ if and only if $x=0$.  Note that the set 
$$\Pi\cap\{(x_1,x_2,x_3,x_4,x_5)\in\rn^5\setminus\{0\}\,|\,P(x_1,x_2,x_3)=\alpha\cdot Q(x_4,x_5)\}$$
is empty.

The map $\hat\Phi$ is invariant under the action of $\rn^+$ on $\rn^5\setminus\{0\}$ so it induces a map $\Phi$ locally defined on the 4-dimensional sphere $S^4$.  Since the natural projection is submersive the fibres $\Phi^{-1}(\{\alpha\})=\hat\Phi^{-1}(\{\alpha\})\cap S^4$ of $\Phi$ are also minimal submanifolds of $S^4$ of codimension two i.e. minimal surfaces.  For a non-zero element $\alpha\in\cn^*$ in the image of $\Phi$ the fibre $\Phi^{-1}(\{\alpha\})$ satisfies the equation $P(x_1,x_2,x_3)=\alpha\cdot Q(x_4,x_5)$.  The {\it minimal surface} $$\Sigma_\alpha^2=\{(x\in S^4\,|\,P(x_1,x_2,x_3)=\alpha\cdot Q(x_4,x_5)\}$$ 
in $S^4$ is closed and bounded and hence {\it compact}. Changing the value of the parameter $\alpha$ will change the shape of the surface $\Sigma^2_\alpha$, so they are not congruent.

\section{Harmonic Morphisms on the Lorentzian Spaces $\rn^n_1$}

By $\rn^n_1=(\rn^n,\ip{\cdot}{\cdot}_L)$ we denote the $n$-dimensional semi-Riemannian  manifold $\rn^n$ equipped with the Lorentzian metric $\ip{\cdot}{\cdot}_L$ given by 
$$\ip{x}{y}_L=-x_n\,y_n+\sum_{k=1}^{n-1}x_k\,y_k.$$
Let $\phi^*,\psi^*:\rn^n_1\to\cn$ be two $C^2$ complex-valued functions, then the corresponding tension field $\tau$ and conformality operator $\kappa$ satisfy
$$\tau(\phi^*)=-\frac{\partial^2\phi^*}{\partial x_n^2}+\sum_{k=1}^{n-1}\frac{\partial^2\phi^*}{\partial x_k^2},$$
$$\kappa(\phi^*,\psi^*)=-\frac{\partial\phi^*}{\partial x_n}\frac{\partial\psi^*}{\partial x_n}+\sum_{k=1}^{n-1}\frac{\partial\phi^*}{\partial x_k}\frac{\partial\psi^*}{\partial x_k}.$$
Let $\phi,\psi:\rn^n\to\cn$ be two complex-valued functions on the Euclidean $\rn^n$ and define their "dual" functions $\phi^*,\psi^*:\rn^n_1\to\cn$ on $\rn^n_1$ by the formulae
$$\phi^*:(x_1,\dots,x_{n-1},x_n)\mapsto\phi(x_1,\dots,x_{n-1},i\cdot x_n),$$
$$\psi^*:(x_1,\dots,x_{n-1},x_n)\mapsto\psi(x_1,\dots,x_{n-1},i\cdot x_n).$$
Then it is clear that $\tau(\phi)=0$ and $\kappa(\phi,\psi)=0$ on the Riemannian $\rn^n$ if and only if $\tau(\phi^*)=0$ and $\kappa(\phi^*,\psi^*)=0$ on the semi-Riemannian $\rn^n_1$.  This duality was first discovered in \cite{Gud-7} and later introduced in the much more general context of Riemannian symmetric spaces in \cite{Gud-Sve-1}.

\section{Harmonic Morphisms on the Hyperbolic Spaces $H^{2n+2}$}

Let $U^{2n+3}$ be the open subset of $\rn^{2n+3}_1$ satisfying 
$$U^{2n+3}=\{x\in\rn^{2n+3}_1\,|\,\ip{x}{x}_L<0\ \ \text{and}\ \ x_{2n+3}>0\}.$$
Further let $H^{2n+2}$ denote the standard model of the hyperbolic space in $\rn^{2n+3}_1$ embedded as follows
$$H^{2n+2}=\{x\in\rn^{2n+3}_1\,|\,\ip{x}{x}_L=-1\ \ \text{and}\ \ x_{2n+3}>0\}.$$
Then we have the radial projection $\pi^*:U^{2n+3}\to H^{2n+2}$ given by $$\pi^*:x\mapsto x/\sqrt{-\ip{x}{x}_L}.$$  Then $\pi^*$ is a well-known submersive harmonic morphism from the semi-Riemannian $U^{2n+3}$ to the Riemannian $H^{2n+2}$.  For this see Lemma 4.1 of \cite{Gud-7}.

The next useful result is a special case of Proposition 1 in \cite{Gud-8}.

\begin{proposition}\label{proposition-invariance-2}
Let $\pi:(M,g)\to(N,h)$ be a non-constant harmonic morphism from a semi-Riemannian manifold $(M,g)$ to a Riemannian $(N,h)$. Further let $\phi:(N,h)\to\cn$ be a complex-valued function and $\hat \phi:(M,g)\to\cn$ be the composition $\hat\phi=\phi\circ\pi$.  Then $\phi$ is a harmonic morphism if and only if $\hat\phi$ is a harmonic morphism.
\end{proposition}

\begin{example}
Let $P_d,Q_d:\cn^n\to\cn$ be the polynomials defined in Example \ref{example-PQ} and $p:\rn^3\to\cn$ as in Example \ref{example-p} with 
$$a_{1}^2 + b_{1}^2 + b_{2}^2=0,\ a_{1}a_{2} = -b_{2}^2\ \ \text{and}\ \ 
a_{1}b_{3} = b_{1} b_{2}.$$	
Set $$z_1=x_4+ix_5,\ \dots\ z_{n-1}=x_{2n}+ix_{2n+1},\ z_{n}=x_{2n+2}-x_{2n+3}$$
and define the map $\hat\Phi^*:V\to\cn$ on the appropriate open subset $V$ of $U^{2n+3}\subset \rn^{2n+3}_1$ by 
$$\hat\Phi^*_d(x_1,\dots,x_n)=\frac{p(x_1,x_2,x_3)^{d/2}+P_d(z_1,\dots,z_n)}{Q_d(z_1,\dots,z_n)}.$$
According to our observations in Section \ref{section-hm-S2n+2} and the above mentioned duality principle all the maps $\hat\Phi_d^*$ are harmonic morphisms.  It should be noted that if $d$ is even and $Q_d^*=z_n^d=(x_{2n+2}-x_{2n+3})^d$ then $\Z(Q_d^*)$ is empty so $\hat\Phi_d^*$ is globally defined on $U^{2n+3}$.
\end{example}

\section{Minimal Surfaces in the Hyperbolic Space $H^{4}$}\label{section-minimal-H4}

Let $a_1,a_2,b_1,b_2,b_3\in\cn$ be complex numbers such that 
$$a_{1}^2 + b_{1}^2 + b_{2}^2=0,\ a_{1}a_{2} = -b_{2}^2,\ a_{1}b_{3} = b_{1} b_{2}\ \ \text{and}\ \ a_1b_1b_2\neq 0.$$
Then define the polynomial functions $P,Q:\rn^5\setminus\{0\}$ by 
$$P(x)=a_{1}(x_1^2-x_2^2)+a_{2}(x_2^2-x_3^2)+b_{1}(x_1x_2)+b_{2}(x_1x_3)+b_{3}(x_2x_3)$$
and 
$$Q(z)=(x_4- x_5)^2.$$
According to our investigations above we know that the quotient $\hat\Phi^*=P/Q$ is a harmonic morphism globally defined on the open subset $U^5$ of $\rn^5_1$.  The complex gradient $\grad(\hat\Phi^*)$
is given by 
$$\grad(\hat\Phi^*)=(
\frac{\partial\hat\Phi^*}{\partial x_1},
\frac{\partial\hat\Phi^*}{\partial x_2},
\frac{\partial\hat\Phi^*}{\partial x_3},
\frac{\partial\hat\Phi^*}{\partial x_4},
\frac{\partial\hat\Phi^*}{\partial x_5}),$$
where 
{\small 
$$\frac{\partial\hat\Phi^*}{\partial x_1}=
\frac{2a_1x_1 + b_1x_2 + b_2x_3}{(x_4 - x_5)^2}
$$
$$\frac{\partial\hat\Phi^*}{\partial x_2}=
\frac{b_1(a_1x_1+2b_1x_2 + b_2x_3)}{a_1(x_4 - x_5)^2}
$$
$$\frac{\partial\hat\Phi^*}{\partial x_3}=
\frac{b_2(a_1x_1 + b_1x_2 + 2b_2x_3)}{a_1(x_4 - x_5)^2}
$$
$$\frac{\partial\hat\Phi^*}{\partial x_4}=
\frac{-2a_1^2(x_1^2 -x_2^2) + 2b_2^2(x_2^2 - x_3^2) - 2a_1b_1x_1x_2 - 2a_1b_2x_1x_3 - 2b_1b_2x_2x_3}{a_1(x_4 - x_5)^3}
$$
$$\frac{\partial\hat\Phi^*}{\partial x_5}=-\frac{\partial\hat\Phi^*}{\partial x_4}.
$$}
Employing the same arguments as in Section \ref{section-minimal-S4} we see that the only solutions to the equation $\grad(\hat\Phi^*)=0$ are the elements of the punctured $(x_4,x_5)$-plane i.e. $$\Pi=\{(0,0,0,x_4,x_5)\,|\,x_4^2+x_5^2\neq 0\}.$$
For a point $x\in\Pi$ we have $\hat\Phi^*(x)=0$.  This shows that the only critical value for $\hat\Phi^*$ is $0\in\cn$.  This means that for any non-zero complex number $\alpha\in\cn^*$ in the image of $\hat\Phi^*$ the fibre $(\hat\Phi^*)^{-1}(\{\alpha\})$ is a regular minimal submanifold of $\rn^5$ of codimension two.

The globally defined map $\hat\Phi^*:U^5\to H^4$ is invariant under the action of $\rn^+$ on $U^5$ so it induces a map $\Phi^*$ globally defined on the 4-dimensional hyperbolic space $H^4$.  Since the natural projection is submersive the fibres $(\Phi^*)^{-1}(\{\alpha\})=(\hat\Phi^*)^{-1}(\{\alpha\})\cap H^4$ of $\Phi^*$ are also minimal submanifolds of $H^4$ of codimension two i.e. {\it minimal surfaces}.  For a non-zero element $\alpha$ in the image of $\Phi^*$ the fibre $(\Phi^*)^{-1}(\{\alpha\})$ satisfies the equation $P(x)-\alpha\cdot Q(z)=0$ so it is {\it complete}.  Changing the value of the parameter $\alpha$ will change the shape of the surface, so they are not congruent.



\begin{thebibliography}{99}

\bibitem{Bai-Eel}
P. Baird, J. Eells,
{\it A conservation law for harmonic maps},
Geometry Symposium Utrecht 1980, Lecture Notes in Mathematics {\bf 894}, Springer (1981), 1-25.


\bibitem{Bai-Woo-book}
P. Baird and J. C. Wood,
{\it Harmonic morphisms between Riemannian manifolds},
London Math. Soc. Monogr. {\bf 29},
Oxford Univ. Press (2003).


\bibitem{Fug-2} B.~Fuglede,
{\it Harmonic morphisms between semi-riemannian manifolds}, Ann.
Acad. Sci. Fennicae {\bf 21} (1996), 31-50.


\bibitem{Gud-4}
S. Gudmundsson,
{\it Harmonic morphisms from complex projective spaces},
Geom. Dedicata {\bf 53} (1994), 155-161.


\bibitem{Gud-7}
S. Gudmundsson,
{\it Minimal submanifolds of hyperbolic spaces via harmonic morphisms},
Geom. Dedicata {\bf 62} (1996), 269-279.

\bibitem{Gud-8}
S.~Gudmundsson, 
{\it On the existence of harmonic morphisms from symmetric spaces of rank one}, Manuscripta Math. {\bf 93} (1997),
421-433.


\bibitem{Gud-bib}
S. Gudmundsson,
{\it The Bibliography of Harmonic Morphisms},
{\tt www.matematik.lu.se/ matematiklu/personal/sigma/harmonic/bibliography.html}


\bibitem{Gud-Sak-1}
S. Gudmundsson, A. Sakovich,
{\it Harmonic morphisms from the classical compact semisimple Lie groups},
Ann. Global Anal. Geom. {\bf 33} (2008), 343-356.

\bibitem{Gud-Sve-1}
S.~Gudmundsson and M.~Svensson,
{\it Harmonic morphisms from the Grassmannians and their non-compact duals},
Ann. Global Anal. Geom. {\bf 30} (2006), 313-333.


\bibitem{Sve-1}
M. Svensson, 
{\it Harmonic morphisms from even-dimensional hyperbolic spaces},
Math. Scand. {\bf 92} (2003), 246–260.


\end{thebibliography}
\end{document}